\documentclass[12pt]{amsart}
\usepackage{amsfonts}
\usepackage{ifthen}
\usepackage{verbatim}
\usepackage{amsthm}
\usepackage{amsmath}
\usepackage{amssymb}
\usepackage{graphicx}
\usepackage{subfigure}
\usepackage{amscd,amssymb,amsthm}
\usepackage{color,xcolor}

\newcounter{minutes}
\divide\time by 60
\newcounter{hours}
\multiply\time by 60 \addtocounter{minutes}{-\time}

\title{ Simpson's Rule revisited}

\author{Slavko Simi\'c}

\address{Mathematical Institute SANU, Kneza Mihaila 36, Belgrade,
Serbia}\email{ssimic@turing.mi.sanu.ac.rs}


\newtheorem{theorem}[equation]{Theorem}

\newtheorem{lemma}[equation]{Lemma}

\newtheorem{remark}[equation]{Remark}

\newcommand{\beq}{\begin{equation}}
\newcommand{\eeq}{\end{equation}}

\numberwithin{equation}{section}

\begin{document}

\def\thefootnote{}
\footnotetext{ \texttt{\tiny File:~\jobname .tex,
          printed: \number\year-\number\month-\number\day,
          \thehours.\ifnum\theminutes<10{0}\fi\theminutes}
} \makeatletter\def\thefootnote{\@arabic\c@footnote}\makeatother

\begin{abstract}
In this article we give some refinements of Simpson's Rule in
cases when it is not applicable in it's classical form i.e., when
the target function is not four times differentiable on a given
interval. Some sharp two-sided inequalities for an extended form
of Simpson's Rule are also proven.
\end{abstract}

\maketitle
\section{Introduction}
\bigskip

We begin with some notions from Classical Analysis which will be
frequently needed in the sequel.

\bigskip

A function $h: I\subset\mathbb R\to \mathbb R$ is said to be
convex on an non-empty interval $I$ if the inequality
\begin{equation}\label{eq1}
h(px+qy)\le ph(x)+qh(y)
\end{equation}
holds for all $x,y\in I$ and all non-negative $p, q; p+q=1$.

If the inequality (\ref{eq1}) reverses, then $h$ is said to be
concave on $I$. \cite{hlp}

The well-known convexity/concavity criteria says that if $h\in
C^{(2)}(I)$ and $h''(x)\gtrless 0, x\in I$, then the function $h$
is convex/concave on $I$. \cite{hlp}

\bigskip

Let $h: I\subset\mathbb R\to \mathbb R$ be a convex function on an
interval $I$ and $a,b\in I$ with $a<b$. Then
\begin{equation}\label{eq2}
h(\frac{a+b}{2})\le\frac{1}{b-a}\int_a^b
h(t)dt\le\frac{h(a)+h(b)}{2}.
\end{equation}

This double inequality is known in the literature as the
Hermite-Hadamard (HH) integral inequality for convex functions.
See, for example, \cite{np} and references therein. There is a
number of refinements and possible generalizations of HH
inequality. Some recent trends can be found in  \cite{wq} and
\cite{se}.

If $h$ is a concave function then both inequalities in (\ref{eq2})
hold in the reversed direction.

\vskip 1cm

Closely connected to the HH inequality is the well-known Simpson's
Rule which is of great importance in numerical integration. It
says that

\begin{lemma} \label{l1} \cite{u} For an integrable function $g$, we have

$$
\int_{x_1}^{x_3} g(t)
dt=\frac{1}{3}h(g_1+4g_2+g_3)-\frac{1}{90}h^5 g^{(4)}(\xi),
 (x_1<\xi<x_3),
$$

where $g_i=g(x_i)$ and $h=:x_2-x_1=x_3-x_2$.
\end{lemma}

Now, by taking $x_1=a, x_2=(a+b)/2, x_3=b$, it follows that
$h=(b-a)/2$ and therefore we obtain another form of Simpson's
Rule:

\begin{equation}\label{eq4}
\frac{g(a)+g(b)}{6}+\frac{2}{3}
g(\frac{a+b}{2})-\frac{1}{b-a}\int_a^b
g(t)dt=\frac{g^{(4)}(\xi)}{2880}(b-a)^4, \ \ a<\xi<b.
\end{equation}

\bigskip

Note that the equation (\ref{eq4}) explicitly supposes that $g\in
C^{(4)}(I)$.

\bigskip

An interesting problem arises if $g\notin C^{(4)}(I)$, i.e., if
$g$ is not a four times continuously differentiable function on
$I$. How to approximate the expression

$$
T_g(a,b):= \frac{g(a)+g(b)}{6}+\frac{2}{3}
g(\frac{a+b}{2})-\frac{1}{b-a}\int_a^b g(t)dt
$$

in this case?

\bigskip

In order to give an answer to this question, we shall consider the
class of functions $g \in C^{(n)}(E)$ which are continuously
differentiable up to $n$-th order on an interval $E:=[a,b]\subset
I$. Since $g^{(n)}(\cdot)$ is a continuous function on a closed
interval, there exist numbers $m_n=m_n(a,b;g):=\min_{t\in
E}g^{(n)}(t)$ and $M_n=M_n(a,b;g):=\max_{t\in E}g^{(n)}(t)$. These
numbers will play an important role in further approximations.

\bigskip

Our task in this article is to demonstrate a method which improves
Simpson's Rule in some characteristic situations.

For example, let $g(\cdot)$ be a twice differentiable function on
$E$. Some preliminary bounds for $T_g(a,b)$ in this case can be
obtained by utilizing Hermite-Hadamard inequality (\ref{eq2}) in a
natural way.

Namely, for a given $g\in C^{(2)}(E)$ define an auxiliary function
$h$ by $h(t):= g(t)-m_2 t^2/2$. Since $h''(t)=g''(t)-m_2\ge 0$, we
see that $h$ is a convex function on $E$. Therefore, applying
Hermite-Hadamard inequality, we get

$$
g(\frac{a+b}{2})-\frac{m_2}{2}(\frac{a+b}{2})^2
$$
$$
\le\frac{1}{b-a}\int_a^b g(t)dt
-\frac{m_2}{2}\frac{b^3-a^3}{3(b-a)}\le
$$
$$
\frac{g(a)+g(b)}{2}-\frac{m_2}{2}\frac{a^2+b^2}{2},
$$

that is,

\begin{equation}\label{eq5}
g(\frac{a+b}{2})+\frac{m_2}{24}(b-a)^2 \le\frac{1}{b-a}\int_a^b
g(t)dt \le \frac{g(a)+g(b)}{2}-\frac{m_2}{12}(b-a)^2.
\end{equation}

\bigskip

On the other hand, taking the auxiliary function $h$ to be
$h(t)=M_2t^2/2-g(t)$, we see that it is also convex on $E$.

Applying Hermite-Hadamard inequality again, we get

\begin{equation}\label{eq6}
\frac{g(a)+g(b)}{2}-\frac{M_2}{12}(b-a)^2\le\frac{1}{b-a}\int_a^b
g(t)dt \le g(\frac{a+b}{2})+\frac{M_2}{24}(b-a)^2.
\end{equation}

\bigskip

Apart from the fact that those inequalities improves HH inequality
in the cases if $g$ is convex $(m_2\ge 0)$ or concave $(M_2\le 0)$
function on $E$, they also lead to an estimation of $T_g(a,b)$, as
follows.

\bigskip

Inequalities $(\ref{eq5})$ and $(\ref{eq6})$ give

$$
\frac{1}{b-a}\int_a^b g(t)dt \le
\frac{g(a)+g(b)}{2}-\frac{m_2}{12}(b-a)^2,
$$
$$
\frac{2}{b-a}\int_a^b g(t)dt \le
2g(\frac{a+b}{2})+\frac{M_2}{12}(b-a)^2.
$$

Hence,

$$
\frac{1}{b-a}\int_a^b g(t)dt-[\frac{g(a)+g(b)}{6}+\frac{2}{3}
g(\frac{a+b}{2})]\le \frac{M_2-m_2}{36}(b-a)^2.
$$

\bigskip

Also, adjusting the left-hand sides of $(\ref{eq5})$ and
$(\ref{eq6})$, we get

$$
2g(\frac{a+b}{2})+\frac{m_2}{12}(b-a)^2 \le\frac{2}{b-a}\int_a^b
g(t)dt,
$$
$$
\frac{g(a)+g(b)}{2}-\frac{M_2}{12}(b-a)^2\le\frac{1}{b-a}\int_a^b
g(t)dt
$$

\bigskip

Therefore,

$$
\frac{1}{b-a}\int_a^b g(t)dt-[\frac{g(a)+g(b)}{6}+\frac{2}{3}
g(\frac{a+b}{2})]\ge -\frac{M_2-m_2}{36}(b-a)^2,
$$

and we finally obtain that

\begin{equation}\label{eq7}
\Bigl|T_g(a,b)\Bigr|=\Big|\frac{g(a)+g(b)}{6}+\frac{2}{3}
g(\frac{a+b}{2})-\frac{1}{b-a}\int_a^b g(t)dt\Big|\le
\frac{1}{36}(M_2-m_2)(b-a)^2,
\end{equation}

whenever $g\in C^{(2)}(E)$.

\bigskip

Another and more efficient method to approximate $T_g(a,b)$ is to
use its integral representations in the cases when $g\in
C^{(1)}(E), g\in C^{(2)}(E)$ or $g\in C^{(3)}(E)$. In this way we
obtain the following estimations:

\bigskip

1. For $g\in C^{(1)}(E), \ \Bigl|T_g(a,b)\Bigr|\le
\frac{5}{72}(b-a)(M_1-m_1)$;

\bigskip

2. For $g\in C^{(2)}(E), \ \Bigl|T_g(a,b)\Bigr|\le
\frac{1}{162}(b-a)^2 (M_2-m_2)$;

\bigskip

3. For $g\in C^{(3)}(E), \
\Bigl|T_g(a,b)\Bigr|\le\frac{1}{1152}(b-a)^3(M_3-m_3)$.

\bigskip

\begin{remark} \ A challenging task and an open problem is to improve the
constants $5/72$ and $1/162$, if possible. We shall prove that
$1/1152$ is the best possible constant in part 3.
\end{remark}

\bigskip

In the sequel we sharply refine Simpson's Rule by assuming that
$f''(\cdot)$ is a convex function on $E$. Then,

$$ 0\le T_g(a,b)\le\frac{(b-a)^2}{162}[\frac{f''(a)+f''(b)}{2}-f''(\frac{a+b}{2})],
$$

( Theorem \ref{thm3}, below).

\bigskip

Finally, applying the method described above, we shall give tight
bounds for an improved form of Simpson's Rule of the fourth order.

\bigskip

\section{Results and Proofs}

\bigskip

We begin with refinements of Simpson's Rule in non-standard cases.

\bigskip

For this cause, the following integral representation of
$T_\phi(a,b)$ in the case $\phi\in C^{(1)}(E)$ is of crucial
value.

\begin{lemma}\label{l0} \ The identity

$$
T_\phi(a,b)=\frac{b-a}{12}\int_0^1 (1-3t)[\phi'(u)-\phi'(v)]dt
$$

holds for any $\phi\in C^{(1)}(E)$, where
$u:=a\frac{t}{2}+b(1-\frac{t}{2}),
v:=b\frac{t}{2}+a(1-\frac{t}{2})$.
\end{lemma}

\bigskip

\begin{proof} \ In the well-known formula

$$
\int_0^1UdV=UV|_0^1-\int_0^1VdU,
$$

putting

$$
U=1-3t, \ dV=(\phi'(u)-\phi'(v))dt,
$$

we get

$$
dU=-3dt, \ V=-\frac{2}{b-a}(\phi(u)+\phi(v)).
$$

\bigskip

Also,

$$
\int_0^1\phi(u)dt=\int_0^1\phi(a\frac{t}{2}+b(1-\frac{t}{2}))dt=\frac{2}{b-a}\int_{\frac{a+b}{2}}^b\phi(t)dt;
$$
$$
\int_0^1\phi(v)dt=\int_0^1\phi(b\frac{t}{2}+a(1-\frac{t}{2}))dt=\frac{2}{b-a}\int_a^{\frac{a+b}{2}}\phi(t)dt.
$$

\bigskip

Therefore,

$$
\frac{b-a}{12}\int_0^1
(1-3t)[\phi'(u)-\phi'(v)]dt=\frac{b-a}{12}[\frac{2}{b-a}(3t-1)(\phi(u)+\phi(v))|_0^1-\frac{6}{b-a}\int_0^1(\phi(u)+\phi(v))dt]
$$
$$
=\frac{2}{3}\phi(\frac{a+b}{2})+\frac{\phi(a)+\phi(b)}{6}-\frac{1}{2}[\frac{2}{b-a}(\int_a^{\frac{a+b}{2}}\phi(t)dt
+\int_{\frac{a+b}{2}}^b\phi(t)dt)]
$$
$$
=\frac{\phi(a)+\phi(b)}{6}+\frac{2}{3}
\phi(\frac{a+b}{2})-\frac{1}{b-a}\int_a^b \phi(t)dt=T_{\phi}(a,b).
$$

\end{proof}

\bigskip

Our first contribution is the following

\begin{theorem}\label{thm0} \ For any $\phi\in C^{(1)}(E)$, we
have

$$
\Big|\frac{\phi(a)+\phi(b)}{6}+\frac{2}{3}
\phi(\frac{a+b}{2})-\frac{1}{b-a}\int_a^b \phi(t)dt\Big|\le
\frac{5}{72}(M_1-m_1)(b-a),
$$

where $m_1:=\min_{t\in E}\phi'(t)$ and $M_1:=\max_{t\in
E}\phi'(t)$.
\end{theorem}

\bigskip

\begin{proof} \ By the above Lemma, we obtain

\bigskip

$$
\Bigl|T_{\phi}(a,b)\Bigr|\le \frac{b-a}{12}\int_0^1
|1-3t||\phi'(u)-\phi'(v)|dt
$$
$$
\le \frac{b-a}{12}(M_1-m_1)\int_0^1
|1-3t|dt=\frac{5}{72}(M_1-m_1)(b-a).
$$

\end{proof}

\bigskip

\begin{theorem}\label{thm1} \ For any $\phi\in C^{(2)}(E)$, we
have

$$
\Big|\frac{\phi(a)+\phi(b)}{6}+\frac{2}{3}
\phi(\frac{a+b}{2})-\frac{1}{b-a}\int_a^b \phi(t)dt\Big|\le
\frac{1}{162}(M_2-m_2)(b-a)^2,
$$

where $m_2:=\min_{t\in E}\phi''(t)$ and $M_1:=\max_{t\in
E}\phi''(t)$.

\end{theorem}

\bigskip

\begin{proof} \ In this case, the integral representation of
$T_\phi(a,b)$ has the following form.

\bigskip

\begin{lemma}\label{l2} \ The identity

$$
T_\phi(a,b)=\frac{(b-a)^2}{48}\int_0^1
t(2-3t)[\phi''(u)+\phi''(v)]dt
$$

holds for any $\phi\in C^{(2)}(E)$, where $u$ and $v$ are defined
as in Lemma \ref{l0}.

\end{lemma}

\begin{proof} \ Indeed, applying partial integration on the
assertion from Lemma \ref{l0}, we get

$$
T_\phi(a,b)=\frac{b-a}{12}\int_0^1 (1-3t)[\phi'(u)-\phi'(v)]dt
$$
$$
= \frac{b-a}{12}[(t-3t^2/2)(\phi'(u)-\phi'(v))_0^1
+\frac{b-a}{2}\int_0^1 (t-3t^2/2)[\phi''(u)+\phi''(v)]dt]
$$
$$
=\frac{(b-a)^2}{48}\int_0^1 t(2-3t)[\phi''(u)+\phi''(v)]dt.
$$

\end{proof}

\bigskip

Hence,

$$
\int_0^1 t(2-3t)[\phi''(u)+\phi''(v)]dt=\int_0^{2/3}
t(2-3t)[\phi''(u)+\phi''(v)]dt-\int_{2/3}^1
t(3t-2)[\phi''(u)+\phi''(v)]dt
$$
$$
\le 2M_2\int_0^{2/3} t(2-3t)dt-2m_2\int_{2/3}^1
t(3t-2)dt=\frac{8}{27}(M_2-m_2),
$$

since $\int_0^{2/3} t(2-3t)dt=\int_{2/3}^1
t(3t-2)dt=\frac{4}{27}$.

\bigskip

Analogously,

$$
\int_0^1 t(2-3t)[\phi''(u)+\phi''(v)]dt=\int_0^{2/3}
t(2-3t)[\phi''(u)+\phi''(v)]dt-\int_{2/3}^1
t(3t-2)[\phi''(u)+\phi''(v)]dt
$$
$$
\ge 2m_2\int_0^{2/3} t(2-3t)dt-2M_2\int_{2/3}^1
t(3t-2)dt=-\frac{8}{27}(M_2-m_2),
$$

and we get the desired result.

\end{proof}

\bigskip

A challenging task is to determine the best possible constant $A$
such that the relation

$$
\Bigl|T_\phi(a,b)\Bigr|\le A(M_2-m_2)(b-a)^2,
$$

holds for any $\phi\in C^{(2)}(E)$.

\bigskip

Note that the function $\phi(\cdot)$, defined on $E=[-x,x]$ by

$$
\phi(t)=\begin{cases}t^3/6 &, t\ge 0;\\
                     -t^3/6 &, t\le 0,

\end{cases}
$$

gives $A\ge 1/288$. Hence $A\in [1/288, 1/162]$.

\bigskip

Another important result concerns the functions which are only
3-times differentiable on $E$.

\bigskip

\begin{theorem}\label{thm2} \ For any $\phi\in C^{(3)}(E)$ we
have

$$
\Big|\frac{\phi(a)+\phi(b)}{6}+\frac{2}{3}
\phi(\frac{a+b}{2})-\frac{1}{b-a}\int_a^b \phi(t)dt\Big|\le
\frac{1}{1152}(M_3-m_3)(b-a)^3.
$$

\bigskip

The constant $C=1/1152$ is best possible.

\end{theorem}

\begin{proof} \ For this case we need a new integral
representation of $T_\phi(a,b)$.

\bigskip

\begin{lemma}\label{l3} \ If $\phi\in C^{(3)}(E)$, then

$$
T_\phi(a,b)=\frac{(b-a)^3}{96}\int_0^1
t^2(1-t)[\phi'''(u)-\phi'''(v)]dt.
$$
\end{lemma}

Indeed, by partial integration we get

$$
\int_0^1 t(2-3t)[\phi''(u)+\phi''(v)]dt=
t^2(1-t)[\phi''(u)+\phi''(v)]|_0^1- \int_0^1
t^2(1-t)\frac{d}{dt}[\phi''(u)+\phi''(v)]
$$
$$
=\frac{b-a}{2}\int_0^1 t^2(1-t)[\phi'''(u)-\phi'''(v)]dt,
$$

and, by Lemma \ref{l2}, the proof follows.

\bigskip

Therefore,

$$
|T_\phi(a,b)|\le \frac{(b-a)^3}{96}\int_0^1
t^2(1-t)|\phi'''(u)-\phi'''(v)|dt
$$
$$
\le \frac{(b-a)^3}{96}(M_3-m_3)\int_0^1 t^2(1-t)dt=
\frac{1}{1152}(M_3-m_3)(b-a)^3.
$$

\bigskip

To prove that the constant $C=1/1152$ is best possible, we
consider the function $d(\cdot)$ defined as:

$$
d(x)=\begin{cases} -x^3/6-x/3,& x\le -1; \\
                    x^4/24+x^2/4-x/6+1/24, & -1\le x\le 1; \\
                    x^3/6, & x\ge 1.
                    \end{cases}
                    $$

 It is easy to confirm that this function is 3-times continuously differentiable
 on the real line.

 Applying the form of Simpson's Rule for $x\in [-a,a], a>1$, we
 obtain

 \begin{equation}\label{eq12}
 \Bigl|\frac{d(-a)+d(a)}{6}+\frac{2}{3}d(0)-\frac{1}{2a}\int_{-a}^a
 d(x)dx \Bigr|\le 8C a^3(M_3-m_3).
 \end{equation}

 Since,

$$
d'''(x)=\begin{cases} -1,& x\le -1; \\
                    x, & -1\le x\le 1; \\
                    1, & x\ge 1,
                    \end{cases}
                    $$

we see that $m_3=-1, \ M_3=1$.

\bigskip

Therefore, by $(\ref{eq12})$ we get

$$
C\ge\frac{|a^3/72-a/36+1/36-1/120a|}{16a^3}
$$
$$
=\frac{1}{1152}\Bigl|1-\frac{2}{a^2}+\frac{2}{a^3}-\frac{3}{5a^4}\Bigr|.
$$

Letting $a\to\infty$, we obtain

$$
C\ge \frac{1}{1152},
$$

and the proof is done.

\end{proof}

\vskip1cm

We shall give in the sequel precise estimation of an extended form
of Simpson's rule under a smoothness condition posed on the target
function.

For example, supposing that $\phi''$ is convex on $E$, we obtain a
clarification of the formula (\ref{eq4}).

\bigskip

\begin{theorem}\label{thm3} \ For a $\phi\in C^{(4)}(E)$, let $\phi''(\cdot)$ be convex on $E$. Then

$$
0\le \frac{\phi(a)+\phi(b)}{6}+\frac{2}{3}
\phi(\frac{a+b}{2})-\frac{1}{b-a}\int_a^b \phi(t)dt
$$
$$
\le\frac{(b-a)^2}{162}
[\frac{\phi''(a)+\phi''(b)}{2}-\phi''(\frac{a+b}{2})].
$$
\end{theorem}

\begin{proof} \ The left-hand side inequality follows from the convexity of $\phi''$ and (\ref{eq4}). For the second inequality we need
an interesting assertion from Convexity Theory \cite{s}.

\begin{lemma}\label{l4} \ Let $h$ be a convex function on $E$
and for some $u,v\in E$, $u+v=a+b$.

\bigskip

Then,

$$
2h(\frac{a+b}{2})\le h(u)+h(v)\le h(a)+h(b).
$$
\end{lemma}

\bigskip

\begin{proof} \ For the left-hand side we have

$$
h(u)+h(v)\ge 2h(\frac{u+v}{2})=2h(\frac{a+b}{2}).
$$

For the right-hand side let $u=pa+qb; p,q\ge 0, p+q=1$. Then
$v=qa+pb$ and, applying (\ref{eq1}), we get

$$
h(u)+h(v)=h(pa+qb)+h(qa+pb)\le
(ph(a)+qh(b))+(qh(a)+ph(b))=h(a)+h(b),
$$

as desired.
\end{proof}

\bigskip

Now, Lemma \ref{l2} gives

$$
T_\phi(a,b)=\frac{(b-a)^2}{48}\int_0^1
t(2-3t)[\phi''(u)+\phi''(v)]dt
$$
$$
=\frac{(b-a)^2}{48}\Bigl(\int_0^{2/3}
t(2-3t)[\phi''(u)+\phi''(v)]dt-\int_{2/3}^1
t(3t-2)[\phi''(u)+\phi''(v)]dt\Bigr).
$$

Since $u,v\in [a,b]$ and $u+v=a+b$, applying Lemma \ref{l4} to
both integrals separately, we obtain

$$
T_\phi(a,b)\le
\frac{(b-a)^2}{48}\Bigl((\phi''(a)+\phi''(b))\int_0^{2/3}
t(2-3t)dt-2\phi''(\frac{a+b}{2})\int_{2/3}^1 t(3t-2)dt\Bigr)
$$
$$
=\frac{(b-a)^2}{324} [\phi''(a)+\phi''(b)-2\phi''(\frac{a+b}{2})],
$$

since $\int_0^{2/3} t(2-3t)dt=\int_{2/3}^1
t(3t-2)dt=\frac{4}{27}$.

\end{proof}

\bigskip

Applying the method which was demonstrated in Introduction, we
obtain an improved form of Simpson's Rule.

\bigskip

\begin{theorem}\label{thm4} \ For any $\psi\in C^{(4)}(E)$, we have

$$
\Bigl|T_\psi(a,b)-\frac{(b-a)^2}{360}[\frac{\psi''(a)+\psi''(b)}{2}-\psi''(\frac{a+b}{2})]\Bigr|
\le\frac{11}{57600}(M_4-m_4)(b-a)^4,
$$

\bigskip

with $m_4=\min_{t\in E}\psi^{(4)}(t),M_4=\max_{t\in
E}\psi^{(4)}(t)$.
\end{theorem}

\begin{proof} \ Take that $\phi(t)=\psi(t)-m_4t^4/24, t\in E$. Since
$\phi^{(4)}(t)=\psi^{(4)}(t)-m_4\ge 0$, we conclude that
$\phi''(\cdot)$ is a convex function on $E$.

Therefore, applying Theorem \ref{thm3} along with the identity

$$
 \frac{a^4+b^4}{6}+\frac{2}{3}
(\frac{a+b}{2})^4-\frac{b^5-a^5}{5(b-a)}=\frac{(b-a)^4}{120},
$$

which follows from (\ref{eq4}), we get

\begin{equation}\label{eq8}
m_4\frac{(b-a)^4}{2880}\le
T_\psi(a,b)\le\frac{(b-a)^2}{162}[\frac{\psi''(a)+\psi''(b)}{2}-\psi''(\frac{a+b}{2})]
-\frac{11}{9}m_4\frac{(b-a)^4}{2880}.
\end{equation}

\bigskip

Consequently, taking $\phi(t)=M_4t^4/24-\psi(t), t\in E$ we have
that $\phi''(\cdot)$ is a convex function on $E$. Hence, using
Theorem \ref{thm3} again, we obtain

\begin{equation}\label{eq9}
\frac{(b-a)^2}{162}[\frac{\psi''(a)+\psi''(b)}{2}-\psi''(\frac{a+b}{2})]
-\frac{11}{9}M_4\frac{(b-a)^4}{2880}\le T_\psi(a,b)\le
M_4\frac{(b-a)^4}{2880}.
\end{equation}

\bigskip

Now, $(\ref{eq8})$ and $(\ref{eq9})$ give

$$
\frac{9}{20}T_\psi(a,b)\le
\frac{(b-a)^2}{360}[\frac{\psi''(a)+\psi''(b)}{2}-\psi''(\frac{a+b}{2})]
-\frac{11}{20}m_4\frac{(b-a)^4}{2880},
$$

and

$$
\frac{11}{20}T_\psi(a,b)\le\frac{11}{20}M_4\frac{(b-a)^4}{2880}.
$$

\bigskip

Adding these inequalities, we obtain

$$
T_\psi(a,b)-\frac{(b-a)^2}{360}[\frac{\psi''(a)+\psi''(b)}{2}-\psi''(\frac{a+b}{2})]
\le\frac{11}{57600}(M_4-m_4)(b-a)^4.
$$

\bigskip

Analogously, adjusting the left-hand sides of $(\ref{eq8})$ and
$(\ref{eq9})$, we get

$$
\frac{9}{20}T_\psi(a,b)\ge
\frac{(b-a)^2}{360}[\frac{\psi''(a)+\psi''(b)}{2}-\psi''(\frac{a+b}{2})]
-\frac{11}{20}M_4\frac{(b-a)^4}{2880},
$$

and

$$
\frac{11}{20}T_\psi(a,b)\ge\frac{11}{20}m_4\frac{(b-a)^4}{2880}.
$$

\bigskip

Therefore,

$$
T_\psi(a,b)-\frac{(b-a)^2}{360}[\frac{\psi''(a)+\psi''(b)}{2}-\psi''(\frac{a+b}{2})]
\ge-\frac{11}{57600}(M_4-m_4)(b-a)^4,
$$

and the proof is done.

\end{proof}

\bigskip

\section{Applications}

\bigskip

As an illustration of results given in this article, we shall
prove the next assertions.

\begin{theorem}\label{thm5} \ For $x>y>0$, we have

$$
\frac{2}{3}-\frac{\coth x-\coth
y}{x-y}-\frac{16}{243}(x^2+xy+y^2)\le\frac{1}{x-y}\int_y^x\frac{\coth
t}{t}dt\le \frac{2}{3}-\frac{\coth x-\coth y}{x-y}.
$$
\end{theorem}

\begin{proof} \ In Theorem \ref{thm3} take $\phi(u)=\phi''(u)=\cosh u, \ u\in
[-2t,2t]$.

We obtain,

\begin{equation}\label{eq10}
\frac{\cosh 2t}{3}+\frac{2}{3}-\frac{8}{81}t^2(\cosh
2t-1)\le\frac{\sinh 2t}{2t}\le \frac{\cosh 2t}{3}+\frac{2}{3}.
\end{equation}

\bigskip

Since $\sinh 2t=2\sinh t\cosh t$ and $\cosh 2t=1+2\sinh^2t$,
dividing both sides of $(\ref{eq10})$ by $\sinh^2 t$, we get

\begin{equation}\label{eq11}
\coth^2t-\frac{1}{3}-\frac{16}{81}t^2\le\frac{\coth
t}{t}\le\coth^2t-\frac{1}{3}.
\end{equation}

Integrating $(\ref{eq11})$ over $t\in [y,x], \ y>0$, the desired
result follows.
\end{proof}

\bigskip

In this case, Theorem \ref{thm4} gives

\begin{theorem}\label{thm6} \ For $x>y>0$, we have

$$
\Bigl|\frac{2}{3}-\frac{\coth x-\coth
y}{x-y}-\frac{4}{45}(x^2+xy+y^2)-\frac{1}{x-y}\int_y^x\frac{\coth
t}{t}dt\Bigr|\le \frac{22}{1125}\frac{x^5-y^5}{x-y}.
$$

\end{theorem}

Proof is left to the reader.

\bigskip

\begin{remark} \ In the same way it is possible to approximate
integrals of the form

$$
\int_a^b\frac{g(t)}{t\tanh t},
$$

where $g(\cdot)$ is a non-negative function on $E$.
\end{remark}

\vskip 1cm

\section{Conclusion}

\bigskip

The results of this paper are of purely theoretical nature.
Namely, we considered the cases when the classical Simpson's Rule
is not applicable, although they are rare in practice. An open
problem of determining best possible constants in Theorems
\ref{thm0} and \ref{thm1} and our solution in Theorem \ref{thm2}
are of the same kind. Comparison of the classical form
$T_{\phi}(a,b)$ given in $(\ref{eq4})$ and new form
$T_{\phi}'(a,b)$ from Theorem \ref{thm4} clearly shows that the
later is much more precise.

For instance,
$$
T_{x^4}(a,b)=(b-a)^4/120, \ T_{x^4}'(a,b)=0; \
T_{x^5}(a,b)=(a+b)(b-a)^4/48, \ T_{x^5}'(a,b)=0.
$$

Further analysis and the composite form of the new Simpson's Rule
is left to the interested reader.

\end{document}